\documentclass[reqno]{amsart}

\newtheorem{thm}{Theorem}[section]

\newtheorem{prop}[thm]{Proposition}
\newtheorem{conj}[thm]{Conjecture}

\theoremstyle{definition}
\newtheorem{defn}[thm]{Definition}

\newtheorem{ques}[thm]{Question}

\def\C{\mathbb{C}}
\def\Z{\mathbb{Z}}

\def\Q{\mathbb{Q}}
\def\R{\mathbb{R}}

\def\CA{{\mathcal A}}

\def\CL{{\mathcal L}}
\def\CM{{\mathcal M}}

\numberwithin{equation}{section}

\title[Floer homology]
{Floer homology in symplectic geometry and \\in mirror symmetry}

\author[Yong-Geun Oh \& Kenji Fukaya]{Yong-Geun Oh \& Kenji Fukaya}
\thanks{Oh thanks A. Weinstein and late A. Floer for putting
everlasting marks on his mathematics. Both authors
thank H. Ohta and K. Ono for a fruitful collaboration on the
Lagrangian intersection Floer theory
which some part of this survey is based on.}

\address{Department of Mathematics, University
of Wisconsin, WI 53706, USA \& Korea Institute for Advanced Study,
Seoul, Korea; oh@math.wisc.edu}
\address{Department of Mathematics,
Kyoto University, Kitashirakawa, Kyoto, Japan; fukaya@math.kyoto-u.ac.jp}

\begin{document}

\begin{abstract}
In this article, the authors review what the Floer homology is and
what it does in symplectic geometry both in the closed string and in
the open string context. In the first case, the authors will
explain how the chain level Floer theory leads to the $C^0$ symplectic
invariants of Hamiltonian flows and to the study of
topological Hamiltonian dynamics. In the second case, the authors
explain how Floer's original construction of Lagrangian
intersection Floer homology is obstructed in general as soon as one leaves the
category of exact Lagrangian submanifolds. They will survey construction,
obstruction and promotion of the Floer complex to the $A_\infty$ category
of symplectic manifolds. Some applications of
this general machinery to the study of the topology of Lagrangian embeddings
in relation to symplectic topology and to mirror symmetry are also reviewed.
\end{abstract}


\keywords{
Floer homology, Hamiltonian flows, Lagrangian submanifolds,
$A_\infty$-structure, mirror symmetry}

\maketitle

\section{Prologue}

The Darboux theorem in symplectic geometry
manifests \emph{flexibility} of the group of
symplectic transformations. On the other hand, the following
celebrated theorem of Eliashberg \cite{eliash1}
revealed subtle \emph{rigidity} of symplectic transformations :
\emph{The
subgroup $Symp(M,\omega)$ consisting of symplectomorphisms
is closed in $Diff(M)$ with respect to the $C^0$-topology.}

This demonstrates that the study of symplectic topology is subtle and
interesting. Eliashberg's theorem relies on a version of non-squeezing
theorem as proven by Gromov \cite{gromov}.
Gromov \cite{gromov} uses the machinery of pseudo-holomorphic
curves to prove his theorem. There is also a different proof
by Ekeland and Hofer \cite{ekel-hofer} of
the classical variational approach to Hamiltonian systems.
The interplay between these two facets of symplectic geometry
has been the main locomotive in the development of symplectic
topology since Floer's pioneering work on his `semi-infinite' dimensional
homology theory, now called the \emph{Floer homology theory}.

As in classical mechanics, there are two most
important boundary conditions in relation to Hamilton's equation
$\dot x = X_H(t,x)$ on a general symplectic manifold :
one is the \emph{periodic} boundary condition
$\gamma(0) = \gamma(1)$,
and the other is the \emph{Lagrangian} boundary condition
$\gamma(0) \in L_0, \, \gamma(1) \in L_1$
for a given pair $(L_0,L_1)$ of two Lagrangian submanifolds
: A submanifold $i: L \hookrightarrow (M,\omega)$ is
called Lagrangian if $i^*\omega = 0$ and $\dim L = \frac{1}{2}\dim M$.
The latter replaces the \emph{two-point} boundary condition in
classical mechanics.

In either of the above two boundary conditions, we have a version of
the \emph{least action principle} : a solution of Hamilton's equation
corresponds to a critical point of the action
functional on a suitable path space with the corresponding
boundary condition. For the periodic boundary condition,
we consider the free loop space
$$
\CL M = \{ \gamma : S^1 \to M\}
$$
and for the Lagrangian boundary condition, we consider the space of
paths connecting
$$
\Omega(L_0, L_1) = \{ \gamma : [0,1] \to M \mid \gamma(0) \in L_0,\,
\gamma(1) \in L_1 \}.
$$
Both $\CL M$ and $\Omega(L_0,L_1)$ have countable number of
connected components. For the case of $\CL M$, it has a distinguished component consisting of
the contractible loops. On the other hand, for the case of
$\Omega(L_0, L_1)$ there is no such distinguished component
in general.
\smallskip
\par
\noindent{\bf Daunting Questions.} For a given time dependent Hamiltonian
$H = H(t,x)$ on $(M,\omega)$, does there exist a solution of the Hamilton equation
$\dot x = X_H(t,x)$ with the corresponding boundary conditions? If so,
how many different ones can exist?
\medskip

One crucial tool for the study of these questions is the
least action principle. Another seemingly trivial but crucial
observation is that when $H \equiv 0$
for the closed case and when $L_1 = L_0$ (and $H \equiv 0$)
for the open case,
there are ``many'' solutions given by \emph{constant} solutions.
It turns out that these two ingredients, combined with Gromov's
machinery of pseudo-holomorphic curves, can be utilized
to study each of the above questions, culminating in
 Floer's proof of Arnold's conjecture for the fixed points
\cite{floer:fixed}, and for the intersection
points of $L$ with its Hamiltonian deformation $\phi_H^1(L)$
 \cite{floer:intersect} for the exact case respectively.

We divide the rest of our exposition into two categories, one in the closed
string and the other in the open string context.

\section{Floer homology of Hamiltonian fixed points}

\subsection{Construction}

On a symplectic manifold $(M,\omega)$, for each given time-periodic
Hamiltonian $H$ i.e., $H$ with $H(t,x) = H(t+1,x)$,
there exists an analog $\CA_H$ to the classical action functional defined on a suitable
covering space of the free loop space. To exploit the fact
that in the vacuum, i.e., when $H \equiv 0$ we have many constant solutions
all lying in the distinguished component of the free
loop space $\CL(M)$
$$
\CL_0(M) = \{ \gamma : [0,1] \to M \mid \gamma(0) = \gamma(1), \,
\mbox{$\gamma$ contractible}\},
$$
one studies the contractible periodic solutions and
so the action functional on $\CL_0(M)$.
The covering space, denoted by $\widetilde\CL_0(M)$, is
realized by the set of suitable equivalence classes $[z,w]$ of the pair
$(z,w)$ where $z:S^1 \to M$ is a loop and $w:D^2 \to M$ is a disc
bounding $z$. Then $\CA_H$ is defined by
\begin{equation}\label{eq:actionH}
\CA_H([z,w]) = - \int w^*\omega - \int_0^1 H(t, \gamma(t))\,
dt.
\end{equation}
This reduces to the classical action $\int pdq - H\, dt$ if we
define the canonical symplectic form as $\omega_0 = \sum_j dq^j \wedge dp_j$
on the phase space $\R^{2n}\cong T^*\R^n$.

To do Morse theory, one needs to introduce a metric on
$\Omega(M)$, which is done by introducing an almost complex
structure $J$ that is compatible to $\omega$ (in that the bilinear
form $g_J:= \omega(\cdot, J \cdot)$
defines a Riemannian metric on $M$) and integrating the norm of the
tangent vectors of the loop $\gamma$. To make the Floer theory
a more flexible tool to use, one should allow this $J$ to be time-dependent.

A computation shows that the negative $L^2$-gradient flow equation of
the action functional for the path $u: \R \times S^1 \to M$
is the following \emph{nonlinear} first
order partial differential equation
\begin{equation}\label{eq:dtaudtH}
\frac{\partial u}{\partial \tau} + J\Big(\frac{\partial
u}{\partial t} - X_H(t,u)\Big) = 0.
\end{equation}
The rest points of this gradient flow are the periodic
orbits of $\dot x = X_H(t,x)$. Note that when $H = 0$, this
equation becomes the pseudo-holomorphic equation
\begin{equation}\label{eq:dtaudu0}
\overline\partial_J(u) =\frac{\partial u}{\partial \tau} +
J\frac{\partial u}{\partial t} = 0
\end{equation}
which has many constant solutions.
Following Floer \cite{floer:fixed},
for each given \emph{nondegenerate} $H$, i.e., one whose
time-one map $\phi_H^1$ has the linearization with no eigenvalue
$1$, we consider a vector space $CF(H)$
consisting of \emph{Novikov Floer chains}
\begin{defn} For each formal sum
\begin{equation}\label{eq:beta}
\beta = \sum _{[z,w] \in \text{Crit}\CA_H} a_{[z,w]} [z,w],
\, a_{[z,w]} \in \Q
\end{equation}
we define the \emph{support} of $\beta$ by the set
$$
\operatorname{supp}\beta = \{ [z,w] \in \operatorname{Crit}\CA_H \mid
a_{[z,w]} \neq 0 \, \mbox{in (\ref{eq:beta})} \}.
$$
We call $\beta$ a Novikov Floer chain or (simply a Floer chain)
if it satisfies the condition
$$
\#\{[z,w] \in \operatorname{supp}\beta \mid
\CA_H([z,w]) \geq \lambda \} < \infty
$$
for all $\lambda \in \R$ and define $CF(H)$ to be the set of
Novikov Floer chains.
\end{defn}

$CF(H)$ can be considered either as a $\Q$-vector space or a module over
the Novikov ring $\Lambda_\omega$ of $(M,\omega)$.
Each Floer chain $\beta$ as a $\Q$-chain
can be regarded as the union of ``unstable manifolds'' of the
generators $[z,w]$ of $\beta$, which has a `peak'.
There is the natural Floer boundary map
$\partial = \partial_{(H,J)} : CF(H) \to CF(H)$
i.e., a linear map satisfying $\partial \partial = 0$. The pair $(CF(H),
\partial_{(H,J)})$ is the \emph{Floer complex} and the quotient
$$
HF_*(H,J;M) := \ker
\partial_{(H,J)}/\operatorname{im}\partial_{(H,J)}
$$
the \emph{Floer homology}. By now the general
construction of this Floer homology has been carried out by Fukaya-Ono \cite{fon},
Liu-Tian \cite{liu-tian1}, and Ruan \cite{ruan} in its
complete generality, after the construction
had been previously carried out by Floer \cite{floer:fixed},
Hofer-Salamon \cite{hofer-sal} and by Ono \cite{ono} in
some special cases.

The Floer homology $HF_*(H,J;M)$ also has the ring structure
arising from the pants product, which
becomes the \emph{quantum product} on $H^*(M)$ in ``vacuum'' i.e.,
when $H \equiv 0$.
The module $H^*(M)\otimes \Lambda_\omega$ with this ring structure is
the \emph{quantum cohomology ring} denoted by $QH^*(M)$. We denote by $a \cdot b$ the quantum
product of two quantum cohomology classes $a$ and $b$.

\subsection{Spectral invariants and spectral norm}

Knowing the existence of periodic orbits of a given Hamiltonian
flow, the next task is to organize the collection of the actions
of different periodic orbits and to study their relationships.

We first collect the actions of all possible periodic orbits,
\emph{including their quantum contributions}, and define the \emph{action
spectrum} of $H$ by
\begin{equation}
\mbox{Spec}(H) : = \{ {\CA}_H([z,w]) \in \R \mid [z,w] \in
\widetilde \Omega_0(M), \, d{\CA}_H([z,w]) = 0\}
\end{equation}
i.e., the set of critical values of ${\CA}_H: \widetilde
\CL_0(M) \to \R$. In general this set is a countable subset of
$\R$ on which the (spherical) period group $\Gamma_\omega$
acts. Motivated by classical Morse theory and mini-max theory,
 one would like to consider
a sub-collection of critical values that are \emph{homologically essential}
: each non-trivial cohomology class gives rise to a mini-max
value, which cannot be pushed further down by the gradient flow.
One crucial ingredient in the classical mini-max theory
is a choice of semi-infinite cycles that are linked under the
gradient flow.

Applying this idea in the context of chain level Floer theory, Oh
generalized his previous construction \cite{oh:jdg,oh:cag} to
the non-exact case in \cite{oh:alan,oh:minimax}.
We define the \emph{level} of a Floer chain
$\beta$ by the maximum value
\begin{equation}\label{eq:lambdaH}
\lambda_H(\beta) : = \max_{[z,w]} \{\CA_H([z,w]) \mid [z,w] \in
\operatorname{supp}\beta\}.
\end{equation}
Now for each $a \in QH^k(M)$ and a generic $J$, Oh considers the mini-max values
\begin{equation}\label{eq:rhoHa}
\rho(H,J;a) = \inf_\alpha \{ \lambda_H(\alpha) \mid \alpha \in CF_{n-k}(H),\,
\partial_{(H,J)} \alpha = 0, \, [\alpha] = a^\flat \}
\end{equation}
where $2n = \dim M$ and proves that this number is independent of $J$ \cite{oh:alan}.
The common number denoted by $\rho(H;a)$ is called the
\emph{spectral invariant} associated to the Hamiltonian $H$
relative to the class $a \in QH^*(M)$. The collection of the
values $\rho(H;a)$ extend to arbitrary smooth Hamiltonian function
$H$, \emph{whether $H$ is nondegenerate or not}, and satisfy the
following basic properties.

\begin{thm}[\cite{oh:alan,oh:minimax}] Let $(M,\omega)$ be an
arbitrary closed symplectic manifold. For any given quantum
cohomology class $0 \neq a \in QH^*(M)$, we have a continuous
function denoted by
$\rho =\rho(H; a): C_m^\infty([0,1] \times M) \times (QH^*(M)\setminus\{0\}) \to \R$
which satisfies the following axioms: Let $H, \, F  \in
C_m^\infty([0,1] \times M)$ be smooth Hamiltonian functions and $a
\neq 0 \in QH^*(M)$. Then we have :
\smallskip \par
$1.$ {\bf (Projective invariance)} $\rho(H;\lambda a) =
\rho(H;a)$ for any $0 \neq \lambda \in \Q$.
\par
$2.$ {\bf (Normalization)} For a quantum cohomology class
$a$, we have $\rho(\underline 0;a) = v(a)$ where $\underline 0$ is the
zero function and $v(a)$ is the valuation of $a$ on $QH^*(M)$.
\par
$3.$ {\bf (Symplectic invariance)} $\rho(\eta^*H ;a) = \rho(H ;a)$
for any $\eta\in Symp(M,\omega)$
\par
$4.$ {\bf (Homotopy invariance)} For any $H,\, K$ with $[H]=[K]$,
$\rho(H;a) = \rho(K;a)$.
\par
$5.$ {\bf (Multiplicative triangle inequality)} $\rho(H \#
F; a\cdot b) \leq \rho(H;a) + \rho(F;b) $
\par
$6.$ {\bf ($C^0$-continuity)} $|\rho(H;a) - \rho(F;a)| \leq
\|H - F\|$. In particular, the function $\rho_a:
H \mapsto \rho(H;a)$ is $C^0$-continuous.
\par
$7.$ {\bf (Additive triangle inequality)} $\rho(H; a+b) \leq \max
\{\rho(H;a), \rho(H;b)\}$
\end{thm}
Under the canonical one-one correspondence between (\emph{smooth}) $H$
(satisfying $\int_M H_t = 0$) and
its Hamiltonian path $\phi_H: t \mapsto \phi_H^t$, we denote by
$[H]$ the path-homotopy class of the Hamiltonian path
$\phi_H: [0,1] \to H am(M,\omega) ; \phi_H(t) = \phi_H^t$ with
fixed end points, and by $\widetilde{H am}(M,\omega)$ the set of
$[H]$'s which represents the universal covering space of
$Ham(M,\omega)$.

This theorem generalizes the results on the exact case by Viterbo
\cite{viterbo2}, Oh \cite{oh:jdg,oh:cag} and Schwarz's
\cite{schwarz} to the non-exact case. The axioms 1 and 7 already
hold at the level of \emph{cycles} or for $\lambda_H$, and follow
immediately from its definition. All other axioms are proved in
\cite{oh:alan} except the homotopy invariance for the
\emph{irrational} symplectic manifolds which is proven in
\cite{oh:minimax}. The additive triangle inequality was explicitly
used by Entov and Polterovich in their construction of some
quasi-morphisms on $Ham(M,\omega)$ \cite{entov-pol1}. The axiom of
homotopy invariance implies that $\rho(\cdot;a)$ projects down to
$\widetilde{Ham}(M,\omega)$. It is a consequence of the following
spectrality axiom, which is proved for \emph{any} $H$ on \emph{rational}
$(M,\omega)$ in \cite{oh:alan} and just for \emph{nondegenerate} $H$ on
\emph{irratioanl} $(M,\omega)$ \cite{oh:minimax} :

\medskip
$8$. {\bf (Nondegenerate spectrality)}
\emph{For nondegenerate $H$, the mini-max values $\rho(H;a)$ lie in $\operatorname{Spec}(H)$
i.e., are critical values of $\CA_H$ for all $a \in QH^*(M) \setminus \{0\}$.}
\medskip

The following is still an open problem.

\begin{ques} Let $(M,\omega)$ be a irrational symplectic manifold,
i.e., the period group $\Gamma_\omega = \{\omega(A) \mid A \in
\pi_2(M)\}$ be a dense subgroup of $\R$. Does $\rho(H;a)$ still
lie in $\operatorname{Spec}(H)$ for all $a \neq 0$ for
\emph{degenerate} Hamiltonian $H$?
\end{ques}

It turns out that the invariant $\rho(H;1)$ can be used to construct
a canonical invariant norm on $Ham(M,\omega)$ of the Viterbo type
which is called the \emph{spectral norm}. To describe this construction,
we start by reviewing the definition of the Hofer norm $\|\phi\|$ of
a Hamiltonian diffeomorphism $\phi$.

There are two natural operations on the space of Hamiltonians $H$
: one the \emph{inverse} $H \mapsto \overline H$ where $\overline
H$ is the Hamiltonian generating the inverse flow
$(\phi_H^t)^{-1}$ and the \emph{product} $(H,F) \mapsto H\# F$
where $H\# F$ is the one generating the composition flow
$\phi_H^t \circ \phi_F^t$. Hofer \cite{hofer} introduced an
invariant norm on $Ham(M,\omega)$. Hofer also considered its
$L^{(1,\infty)}$-version $\|\phi\|$ defined by
$$
\|\phi\| = \inf_{H\mapsto \phi} \|H\| \, ;\quad
\|H\| = \int_0^1(\max H_t - \min H_t) \, dt
$$
where $H \mapsto \phi$ stands for $\phi = \phi_H^1$.
We call $\|H\|$ the \emph{$L^{(1,\infty)}$-norm} of $H$ and
$\|\phi\|$ the \emph{$L^{(1,\infty)}$ Hofer norm} of $\phi$.

Using the spectral invariant $\rho(H;1)$, Oh \cite{oh:dmj}
defined a function $\gamma: C^\infty_m([0,1] \times M) \to \R$ by
$$
\gamma(H) = \rho(H;1) + \rho(\overline H;1)
$$
on $C^\infty_m([0,1]\times M)$, whose definition is more
\emph{topological} than $\|H\|$. For example, $\gamma$ canonically
projects down to a function on $\widetilde{Ham}(M,\omega)$ by
the homotopy invariance axiom while $\|H\|$ does not.
Obviously $\gamma(H) = \gamma(\overline H)$.
The inequality $\gamma(H) \leq \|H\|$
was also shown in \cite{oh:jdg,oh:dmj} and
the inequality  $\gamma(H) \geq 0$ follows from
the triangle inequality applied to $a=b=1$ and from the
normalization axiom $\rho(\underline{0};1) = 0$.

Now we define a non-negative function $\gamma: H am(M,\omega) \to
\R_+$ by
$
\gamma(\phi):= \inf_{H \mapsto \phi}\gamma(H).
$
Then the following theorem is proved in \cite{oh:dmj}.

\begin{thm}[\cite{oh:dmj}] Let $(M,\omega)$ be any closed
symplectic manifold. Then $\gamma : H am(M,\omega) \to \R_+$
defines a (non-degenerate) norm on $Ham(M,\omega)$ which satisfies
the following additional properties :
\par
$1.$ $\gamma(\eta^{-1}\phi \eta) =
\gamma(\phi)$ for any symplectic diffeomorphism $\eta$
\par
$2.$ $\gamma(\phi^{-1}) = \gamma(\phi)$, \quad $\gamma(\phi) \leq \|\phi\|$.
\end{thm}
Oh then applied the function $\gamma = \gamma(H)$ to the study of the
geodesic property of Hamiltonian flows \cite{oh:ajm2,oh:dmj}.

As another interesting application of spectral invariants is
a new construction of \emph{quasi-morphisms} on $Ham(M,\omega)$
carried out by Entov and Polterovich \cite{entov-pol1}.
Recall that for a closed $(M,\omega)$, there exists no non-trivial
homomorphism to $\R$ because $Ham(M,\omega)$ is a simple group \cite{banyaga}.
However for a certain class of \emph{semi-simple}
symplectic manifolds, e.g. for
$ (S^2, \omega),\, (S^2 \times S^2, \omega\oplus \omega), \,
(\C P^n, \omega_{FS})$,
Entov and Polterovich \cite{entov-pol1} produced
non-trivial quasi-morphisms, exploiting the spectral invariants
$\rho(e,\cdot)$ corresponding to a certain idempotent element
$e$ of the quantum cohomology ring $QH^*(M)$.

It would be an important problem to unravel what the true meaning of
Gromov's pseudo-holomorphic curves or of the Floer homology in general
is in regard to the underlying symplectic structure.

\subsection{Towards topological Hamiltonian dynamics}

We note that construction of
spectral invariants largely depends on the smoothness (or at least
differentiability) of Hamiltonians $H$ because it involves the study
of Hamilton's equation $\dot x = X_H(t,x)$.
If $H$ is smooth, there is a one-one correspondence between
$H$ and its flow $\phi_H^t$. However this correspondence
breaks down  when $H$ does not have enough regularity,
e.g., if $H$ is only continuous or even $C^1$
\emph{because the fundamental existence and uniqueness
theorems of ODE fail}.

However the final outcome $\rho(H;a)$ still
makes sense for and can be extended to a certain
natural class of $C^0$-functions $H$.
Now a natural questions to ask is

\begin{ques} Can we define the notion of topological
Hamiltonian dynamical systems? If so, what is the dynamical
meaning of the numbers $\rho(H;a)$ when $H$ is just continuous but
not differentiable?
\end{ques}
These questions led to the notions of
\emph{topological Hamiltonian paths} and
\emph{Hamiltonian homeomorphisms} in \cite{oh:hameo1}.
\begin{defn}\label{topohamflow} A continuous path $\lambda: [0,1] \to Homeo(M)$ with
$\lambda(0) = id$ is called a topological Hamiltonian path if
there exists a sequence of smooth Hamiltonians $H_i:[0,1] \times M
\to \R$ such that
\par
$1.$ $H_i$ converges in the $L^{(1,\infty)}$-topology (or Hofer topology)
of Hamiltonians and
\par
$2.$ $\phi_{H_i}^t \to \lambda(t)$ uniformly converges on $[0,1]$.
\end{defn}
We say that the $L^{(1,\infty)}$-limit of any such sequence
$H_i$ is a \emph{Hamiltonian} of the topological Hamiltonian
flow $\lambda$. The following uniqueness result is proved in
\cite{oh:hameo2}.
\begin{thm}[\cite{oh:hameo2}]
Let $\lambda$ be a topological Hamiltonian path. Suppose that there
exist two sequences $H_i$ and $H_i'$ satisfying the conditions in
Definition \ref{topohamflow}. Then their limits
coincide as an $L^{(1,\infty)}$-function.
\end{thm}
The proof of this theorem is a modification of Viterbo's proof
\cite{viterbo3} of a similar uniqueness result for the $C^0$ Hamiltonians,
combined with a structure theorem of topological Hamiltonians
which is also proven in \cite{oh:hameo2}.
An immediate corollary is the following extension of the spectral invariants
to the space of topological Hamiltonian paths.

\begin{defn} Suppose $\lambda$ is a topological Hamiltonian path and
let $H_i$ be the sequence of smooth Hamiltonians that converges in
$L^{(1,\infty)}$-topology and whose associated Hamiltonian paths $\phi_{H_i}$
converges to $\lambda$ uniformly. We define
$$
\rho(\lambda;a) = \lim_{i \to \infty}\rho(H_i;a).
$$
\end{defn}
The uniqueness theorem of topological Hamiltonians and the $L^{(1,\infty)}$
continuity property of $\rho$
imply that this definition is well-defined.

\begin{defn} A homeomorphism $h$ of $M$ is a Hamiltonian homeomorphism
if there exists a sequence of smooth Hamiltonians $H_i:[0,1] \times M
\to \R$ such that
\par
$1.$ $H_i$ converges in the $L^{(1,\infty)}$-topology
of Hamiltonians and
\par
$2.$ the Hamiltonian path $\phi_{H_i}: t \mapsto \phi_{H_i}^t$ uniformly
converges on $[0,1]$ in the $C^0$-topology of $Homeo(M)$, and $\phi_{H_i}^1 \to h$.
\par
We denote by $Hameo(M,\omega)$ the set of such homeomorphisms.
\end{defn}
Motivated by Eliashberg's rigidity theorem, we also define the group
$Sympeo(M,\omega)$ as the subgroup of $Homeo(M)$ consisting
of the $C^0$-limits of symplectic diffeomorphisms.
Then Oh and M\"uller \cite{oh:hameo1} proved the following theorem

\begin{thm}[\cite{oh:hameo1}] $Hameo(M,\omega)$ is a path-connected
normal subgroup of $Sympeo_0(M,\omega)$, the identity component of
$Sympeo(M,\omega)$.
\end{thm}
One can easily derive that $Hameo(M,\omega)$ is a proper subgroup of
$Sympeo_0(M,\omega)$ whenever the so called mass flow homomorphism
\cite{fathi} is non-trivial or there exists a symplectic diffeomorphism
that has no fixed point, e.g., $T^{2n}$ \cite{oh:hameo1}.
In fact, we conjecture that this is always the case.

\begin{conj} The group $Hameo(M,\omega)$ is a proper subgroup of
$Sympeo_0(M,\omega)$ for any closed symplectic manifold $(M,\omega)$.
\end{conj}
A case of particular interest is the case $(M,\omega) = (S^2,\omega)$.
In this case, together with the smoothing result proven in \cite{oh:smoothing},
the affirmative answer to this conjecture would
negatively answer to the following open question in the area preserving
dynamical systems. See \cite{fathi} for the basic theorems on the
measure preserving homeomorphisms in dimension greater than equal to 3.

\begin{ques} Is the identity component of the group of area preserving
homeomorphisms on $S^2$ a simple group?
\end{ques}

\section{Floer theory of Lagrangian intersections}

Floer's original definition \cite{floer:intersect} of the homology
$HF(L_0,L_1)$ of Lagrangian submanifolds
meets many obstacles when one attempts to generalize his definition
beyond the exact cases i.e., the case
$$
L_0 = L, \quad L_1 = \phi(L) \quad \mbox{with }\, \pi_2(M,L) = \{0\}.
$$
In this exposition, we will consider the cases of Lagrangian submanifolds
that are not necessarily exact. In the open string case
of Lagrangian submanifolds, one has to deal with the phenomenon of bubbling-off discs
besides bubbling-off spheres. One crucial difference between
the former and the latter is that the former is a phenomenon of
codimension one while the latter is that of
codimension two. This difference makes the general Lagrangian intersection Floer theory
display very different perspective compared to the Floer theory of Hamiltonian
fixed points. For example, for the intersection case in general, one has to
study the theory \emph{in the chain level}, which forces one to
consider the chain complexes. Then the
meaning of invariance of the resulting objects is much more non-trivial
to define compared to that of Gromov-Witten invariants for which one can work
with in the level of homology.

There is one particular case that Oh  singled out in \cite{oh:cpam}
for which the original version of Floer cohomology is well-defined and invariant just under the
change of almost complex structures and under the Hamiltonian isotopy.
This is the case of \emph{monotone}
Lagrangian submanifolds with
\emph{minimal Maslov number $\Sigma_L \geq 3$} :
\begin{defn}\label{monotone} A Lagrangian submanifold
$L \subset (M,\omega)$ is \emph{monotone} if there exists a constant $\lambda \geq 0$
such that $\omega(A) = \lambda \mu(A)$ for all elements $A \in \pi_2(M,L)$.
The minimal Maslov number is defined by the integer
$$
\Sigma_L = \min \{ \mu(\beta) \mid \beta \in \pi_2(M,L), \, \mu(\beta) > 0 \}.
$$
\end{defn}
We will postpone further discussion on this particular case until later in this survey
but proceed with describing the general story now.

To obtain the maximal possible generalization of Floer's construction,
it is crucial to develop a proper \emph{off-shell} formulation of
the relevant Floer moduli spaces.

\subsection{Off-shell formulation}
\label{subsec:off-shell}

We consider the space of paths
$$
\Omega = \Omega(L_0,L_1) = \{ \ell :[0,1] \to P~|~ \ell(0) \in L_0,
\ell(1)\in L_1 \}.
$$
On this space, we are given the \emph{action one-form} $\alpha$ defined
by
$$
\alpha(\ell)(\xi) = \int_0^1 \omega(\dot \ell(t),\xi(t)) \, dt
$$
for each tangent vector $\xi \in T_\ell\Omega$. From this expression, it follows
that
$$
\operatorname{Zero}(\alpha) = \{\ell_p: [0,1] \to M \mid p \in L_0\cap L_1,
\quad  \ell_p\equiv p \}.
$$
Using the Lagrangian property of
$(L_0,L_1)$, a straightforward calculation shows that this form is
\emph{closed}. Note that $\Omega(L_0,L_1)$ is not connected but has countably many
connected components. We will  work on a particular fixed connected
component of $\Omega(L_0,L_1)$. We pick up a based path
$\ell_0 \in \Omega(L_0,L_1)$ and consider the corresponding
component $\Omega(L_0,L_1;\ell_0)$, and then
define a covering space
$$
\pi:\widetilde \Omega(L_0,L_1;\ell_0) \to \Omega(L_0,L_1;\ell_0)
$$
on which we have a single valued action functional such that
$d\CA = - \pi^*\alpha$.
One can repeat Floer's construction similarly as in
the closed case replacing $\CL_0(M)$ by the chosen component
of the path space $\Omega(L_0,L_1)$.
We refer to \cite{FOOO} for the details of this construction. We then denote
by $\Pi(L_0,L_1;\ell_0)$ the group of deck transformations.
We define the associated Novikov ring $\Lambda (L_0,L_1;\ell_0)$
as a completion of the group ring $\Q[\Pi(L_0,L_1;\ell_0)]$.

\begin{defn} $\Lambda_k (L_0,L_1;\ell_0)$ denotes the set of
all (infinite) sums
$$\sum_{g\in \Pi(L_0,L_1;\ell_0) \atop \mu (g) = k} a_g
[g]$$ such that $a_g \in \Q$ and that for each $C \in \R$, the set
$$
\# \{ g \in \Pi(L_0,L_1;\ell_0) \mid E(g) \leq C, \,\,a_g \not = 0\}
< \infty.
$$
We put $\Lambda (L_0, L_1;\ell_0) = \bigoplus_k \Lambda_k (L_0, L_1;\ell_0)$.
\end{defn}
We call this graded ring the \emph{Novikov ring} of the pair $(L_0,L_1)$
relative to the path $\ell_0$. Note that
this ring depends on $L$ and $\ell_0$. In relation to mirror symmetry, one needs to
consider a family of Lagrangian submanifolds and to use a universal
form of this ring. The following ring was introduced in
\cite{FOOO} which plays an important role in the rigorous
formulation of homological mirror symmetry conjecture.

\begin{defn}[Universal Novikov ring]
We define
\begin{eqnarray}
\Lambda_{nov} & = & \left\{\sum_{i=1}^\infty a_i T^{\lambda_i} ~\Big|~
a_i \in \Q, \, \lambda_i \in \R, \, \lambda_i \leq \lambda_{i+1}, \,
\lim_{i\to \infty}\lambda_i = \infty \right\}\\
\Lambda_{0,nov} & = &\left\{\sum_{i=1}^\infty a_i T^{\lambda_i} \in \Lambda_{nov}
~\Big|~ \lambda_i \geq 0 \right\}.
\end{eqnarray}
\end{defn}
In the above definitions of Novikov rings, one can replace $\Q$ by other
commutative ring with unit, e.g., $\Z$, $\Z_2$ or $\Q[e]$ with
a formal variable $e$.

There is a natural filtration on these rings provided by
the valuation $v : \Lambda_{nov}, \, \Lambda_{0,nov} \to \R$ defined by
\begin{equation}\label{eq:vlambda1}
v\left(\sum_{i=1}^\infty a_i T^{\lambda_i}\right): = \lambda_1.
\end{equation}
This is well-defined by the definition of the Novikov ring and induces a filtration
$F^\lambda\Lambda_{nov}: = v^{-1}([\lambda, \infty))$
on $\Lambda_{nov}$. The function $e^{-v}: \Lambda_{nov} \to \R_+$ also provides a natural
non-Archimedean norm on $\Lambda_{nov}$.
We call the induced topology on $\Lambda_{nov}$
a \emph{non-Archimedean topology}.

We now assume that $L_0$ intersects $L_1$ transversely and form the
$\Q$-vector space $CF(L_0,L_1)$ over the set
$\operatorname{span}_\Q\operatorname{Crit}\CA$  similarly as $CF(H)$.
Now let $p, \, q \in L_0 \cap L_1$. We
denote by $\pi_2(p,q)=\pi_2(p,q;L_0,L_1)$ the set of homotopy
classes of smooth maps $u: [0,1] \times [0,1]  \to M$
relative to the boundary
$$
u(0,t) \equiv p , \quad u(1,t) = q; \quad u(s,0) \in L_0, \quad
u(s,1) \in L_1
$$
and by $[u] \in \pi_2(p,q)$ the homotopy class of $u$ and by $B$ a
general element in $\pi_2(p,q)$.
For given $B \in \pi_2(p,q)$, we denote by ${Map}(p,q;B)$
the set of such $w$'s in class $B$. Each element $B \in \pi_2(p,q)$ induces a map
given by the obvious gluing map $[p,w] \mapsto [q,w \# u]$
for $u \in Map(p,q;B)$. There is also the natural gluing map
$$
\pi_2(p,q) \times \pi_2(q,r) \to \pi_2(p,r)
$$
induced by the concatenation $(u_1, u_2) \mapsto u_1\# u_2$.

\subsection{Floer moduli spaces and Floer operators}

Now for each given $J = \{J_t\}_{0\leq t \leq 1}$
and $B \in \pi_2(p,q)$, we define the moduli space
$\widetilde\CM(p,q;B)$
consisting of finite energy solutions of the Cauchy-Riemann equation
$$
\begin{cases}
\frac{du}{d\tau} + J_t \frac{du}{dt} = 0 \\
u(\tau,0) \in L_0, \quad u(\tau,1) \in L_1, \, \int u^*\omega < \infty
\end{cases}
$$
with the asymptotic condition and the homotopy condition
$$
u(-\infty,\cdot) \equiv p,
\quad u(\infty,\cdot) \equiv q; \quad [u] = B.
$$
We then define $\CM(p,q;B) = \widetilde\CM(p,q;B)/\R$ the quotient by
the $\tau$-translations and a collection of \emph{rational} numbers
$n(p,q;J,B) = \#(\CM(p,q;J,B))$
whenever the expected dimension of $\CM(p,q;B)$ is zero.
Finally we define the Floer `boundary' map $\partial : CF(L_0,L_1;\ell_0) \to
CF(L_0,L_1;\ell_0)$ by the sum
\begin{equation}\label{eq:boundary}
\partial ([p,w]) = \sum_{q \in L_0\cap L_1}\sum_{B \in \pi_2(p,q)}
n(p,q;J,B) [q,w\# B].
\end{equation}

When  a Hamiltonian isotopy $\{L'_s\}_{0 \leq s\leq 1}$ is given
one also considers the non-autonomous version of the Floer equation
$$
\begin{cases}
\frac{du}{d\tau} + J_{t,\rho(\tau)} \frac{du}{dt} = 0 \\
u(\tau,0) \in L, \quad u(\tau,1) \in L'_{\rho(\tau)}
\end{cases}
$$
as done in \cite{oh:cpam} where $\rho:\R \to [0,1]$ is a smooth
function with $\rho(-\infty) = 0, \, \rho(\infty) = 1$
such that $\rho'$ is compactly supported and define
the Floer `chain' map
$$
h: CF^*(L_0,L_0') \to CF^*(L_1,L'_1).
$$
However unlike the closed case or the exact case, many things go wrong when
one asks for the property $\partial\circ \partial = 0$ or
$\partial h + h \partial = 0$ especially over the rational coefficients,
and even when $HF^*(L,\phi_H^1(L))$ is defined, it is not isomorphic to
the classical cohomology $H^*(L)$.

In the next 3 subsections, we explain how to overcome these troubles
and describe the spectral sequence relating $HF^*(L,\phi_H^1(L))$
to $H^*(L)$ when the former is defined. All the results in these subsections are joint
works with H. Ohta and K. Ono that appeared in \cite{FOOO},
unless otherwise said. We refer to Ohta's article \cite{ohta}
for a more detailed survey on the work from \cite{FOOO}.

\subsection{Orientation}
\label{subsec:orientation}

We first recall the following definition from \cite{FOOO}.

\begin{defn} A submanifold $L \subset M$ is called
\emph{relatively spin} if it is orientable and
there exists a class $st \in H^2(M,\Z_2)$ such that
$st|_L = w_2(TL)$ for the Stiefel-Whitney class $w_2(TL)$ of $TL$.
A pair $(L_0,L_1)$ is relatively spin, if there exists a class
$st \in H^2(M,\Z_2)$ satisfying $st|_{L_i} = w_2(TL_i)$ for each
$i = 0, \, 1$.
\end{defn}

We fix such a class $st \in H^2(M,\Z_2)$ and a triangulation of $M$.
Denote by $M^{(k)}$ its $k$-skeleton. There exists a unique
rank 2 vector bundle $V(st)$ on $M^{(3)}$ with
$w^1(V(st)) = 0, \, w^2(V(st)) = st$. Now suppose that $L$ is
relatively spin and $L^{(2)}$ be the 2-skeleton of $L$.
Then $V\oplus TL$ is trivial on the 2-skeleton
of $L$. We define

\begin{defn} We define a $(M,st)$-relative
spin structure of $L$ to be a spin structure of the
restriction of the vector bundle $V\oplus TL$ to $L^{(2)}$.
\end{defn}

The following theorem was proved by de Silva \cite{silva}
and in \cite{FOOO} independently.

\begin{thm} The moduli space of
pseudo-holomorphic discs is orientable, if $L\subset (M,\omega)$
is relatively spin Lagrangian submanifold. Furthermore the choice
of relative spin structure on $L$ canonically determines an orientation on
the moduli space $\CM(L;\beta)$ of holomorphic discs
for all $\beta \in \pi_2(M,L)$.
\end{thm}

For the orientations on the Floer moduli spaces, the following
theorem was proved in \cite{FOOO}.

\begin{thm} Let $J = \{J_t\}_{0 \leq t \leq 1}$ and
suppose that a pair of Lagrangian submanifolds
$(L_0, L_1)$ are $(M,st)$-relatively spin.
Then for any $p, q \in L_0\cap L_1$ and $B \in \pi_2(p,q)$,
the Floer moduli space $\CM(p,q;B)$ is orientable.
Furthermore a choice of relative spin structures for
the pair $(L_0,L_1)$ determines an orientation on $\CM(p,q;B)$.
\end{thm}

One can amplify the orientation to the moduli space of
pseudo-holomorphic polygons $\CM(\CL,\vec p;B)$ where
$\CL = (L_0, L_1, \cdots, L_k)$ and $\vec p = (p_{01},p_{12}, \cdots,
p_{k0})$ with $p_{ij} \in L_i \cap L_j$ and extend the
construction to the setting of  $A_\infty$ category \cite{Ainfold}.
We refer to \cite{MsurvII} for more detailed discussion on this.

\subsection{Obstruction and $A_\infty$ structure}
\label{subsec:obstruction}

Let $(L_0,L_1)$ be a relatively spin pair with $L_0$ intersecting
$L_1$ transversely
and fix a $(M,st)$-relatively spin structure on each $L_i$.
To convey the appearance of obstruction to the boundary
property $\partial \partial = 0$ in a coherent way,
we assume in this survey, for the simplicity, that
all the Floer moduli spaces involved in the construction are
transverse and so the expected dimension is the same as
the actual dimension. For example, this is the case
for monotone Lagrangian submanifolds at least for the Floer
moduli spaces of dimension 0, 1 and 2. However, we would like to emphasize
that we have to use the machinery of \emph{Kuranishi structure}
introduced in \cite{fon} in the level of \emph{chain} to properly treat
the transversality problem for the general case, whose detailed
study we refer to \cite{FOOO}.

We compute $\partial\partial ([p,w])$. According to the definition
(\ref{eq:boundary}) of the map $\partial$, we have the formula for its matrix
coefficients
\begin{equation}\label{eq:B=B1B2}
\langle \partial\partial [p,w], [r,w\# B] \rangle = \sum_{q \in L_0\cap L_1}
\sum_{B = B_1\# B_2 \in \pi_2(p,r)} n(p,q;B_1)n(q,r;B_2)
\end{equation}
where $B_1 \in \pi_2(p,q)$ and $B_2 \in \pi_2(q,r)$. To prove,
$\partial \partial = 0$, one needs to prove
$\langle \partial\partial [p,w], [r,w\#B] \rangle = 0$
for all pairs $[p,w], \, [r,w\#B]$. On the other hand it follows from
definition that each summand $n(p,q;B_1)n(q,r;B_2)$ is nothing but
the number of broken trajectories lying in
$\CM(p,q;B_1) \# \CM(q,r;B_2)$.
The way how Floer \cite{floer:intersect} proved the vanishing of (\ref{eq:B=B1B2})
under the assumption
\begin{equation}\label{eq:pi2MLi}
L_0 = L, \, L_1 = \phi_H^1(L) ; \, \pi_2(M,L_i) = 0
\end{equation}
is to construct a suitable compactification
of the one-dimensional (smooth) moduli space $\CM(p,r;B) = \widetilde \CM(p,r;B)/\R$
in which the broken trajectories of the form $u_1 \# u_2$ comprise \emph{all} the boundary
components of the compactified moduli space. By definition, the expected
dimension of $\CM(p,r;B)$ is one and so compactified moduli space
becomes a compact one-dimensional manifold. Then
$\partial\partial = 0$ follows.

As soon as one goes beyond Floer's case (\ref{eq:pi2MLi}), one must
concern the problems of \emph{a priori energy bound} and
\emph{bubbling-off discs}.
As in the closed case, the Novikov ring
is introduced to solve the problem of energy bounds. On the other hand,
bubbling-off-discs is a new phenomenon which is that of codimension
one and can indeed occur in the boundary of the
compactification of Floer moduli spaces.

To handle the problem of bubbling-off discs, Fukaya-Oh-Ohta-Ono
\cite{FOOO}, associated a structure of filtered $A_\infty$
algebra $(C, \mathfrak m)$ \emph{with non-zero $\mathfrak m_0$-term} in general. The notion
of $A_\infty$ structure was first introduced by Stasheff \cite{stasheff}.
We refer to \cite{GJ} for an exposition close to ours with different
sign conventions.
The above mentioned obstruction is closely related to non-vanishing of
$\mathfrak m_0$ in this $A_\infty$ structure. Description of this
obstruction is now in order.

Let $C$ be a graded $R$-module where $R$ is the coefficient ring.
In our case, $R$ will be $\Lambda_{0,nov}$. We denote by $C[1]$
its suspension defined by $C[1]^k = C^{k+1}$.
We denote by $deg(x)=|x|$ the degree of $x \in C$ before the shift and
$deg'(x)=|x|'$ that after the degree shifting, i.e., $|x|' = |x| - 1$.
Define the {\it bar complex} $B(C[1])$ by
$$
B_k(C[1]) = (C[1])^{k\otimes}, \quad B(C[1]) =
\bigoplus_{k=0}^\infty B_k(C[1]).
$$
Here $B_0(C[1]) = R$ by definition. We provide the degree of
elements of $B(C[1])$ by the rule
\begin{equation}\label{eq:degonBC[1]}
|x_1 \otimes \cdots \otimes x_k|': = \sum_{i=1}^k |x_i|' = \sum_{i =1}^k|x_i|
-k
\end{equation}
where $|\cdot|'$ is the shifted degree. The ring $B(C[1])$
has the structure of {\it graded coalgebra}.

\begin{defn} The structure of (strong) {\it $A_\infty$ algebra} is a
sequence of $R$ module homomorphisms
$$
\mathfrak m_k: B_k(C[1]) \to C[1], \quad k = 1, 2, \cdots,
$$
of degree +1 such that the coderivation
$d = \sum_{k=1}^\infty \widehat{\mathfrak m}_k$
satisfies $d d= 0$, which
is called the \emph{$A_\infty$-relation}.
Here we denote by $\widehat{\mathfrak
m}_k: B(C[1]) \to B(C[1])$ the unique extension of $\mathfrak m_k$
as a coderivation on $B(C[1])$. A \emph{filtered $A_\infty$ algebra}
is an $A_\infty$ algebra with a filtration for which $\frak m_k$ are
continuous with respect to the induce non-Archimedean topology.
\end{defn}

In particular, we have $\mathfrak m_1
\mathfrak m_1 = 0$ and so it defines a complex $(C,\mathfrak m_1)$. We
define the $\mathfrak m_1$-cohomology by
\begin{equation}\label{eq:m1cohom}
H(C,\mathfrak m_1) = \mbox{ker }\mathfrak m_1/\mbox{im }\mathfrak m_1.
\end{equation}
A {\it weak $A_\infty$ algebra} is defined in the same way, except
that it also includes
$$
\mathfrak m_0: R \to B(C[1]).
$$
The first two terms of the $A_\infty$ relation for a weak
$L_\infty$ algebra are given as
\begin{eqnarray}
\mathfrak m_1(\mathfrak m_0(1)) & = & 0 \label{eq:m1m0=0} \\
\mathfrak m_1\mathfrak m_1 (x) + (-1)^{|x|'}\mathfrak m_2(x, \mathfrak m_0(1)) +
\mathfrak m_2(\mathfrak m_0(1), x) & = & 0. \label{eq:m0m1}
\end{eqnarray}
In particular, for the case of weak $A_\infty$ algebras, $\mathfrak
m_1$ will not necessarily satisfy the boundary property, i.e., $\mathfrak m_1\mathfrak m_1
\neq 0$ in general.

The way how a weak $A_\infty$ algebra is attached to a Lagrangian
submanifold $L \subset (M,\omega)$ arises as an $A_\infty$ deformation of
the classical singular cochain complex including the instanton
contributions. In particular, when there is no instanton contribution
as in the case $\pi_2(M,L) = 0$, it will reduce to an $A_\infty$ deformation
of the singular cohomology in the chain level including all possible
higher Massey product. One outstanding circumstances arise in relation to
the \emph{quantization} of rational homotopy theory on the cotangent bundle $T^*N$
of a compact manifold $N$. In this case, the authors proved in \cite{foh}
that the $A_\infty$ sub-category `generated' by such graphs is literally
isomorphic to a certain $A_\infty$ category constructed by
the Morse theory of \emph{graph flows}.

We now describe the basic $A_\infty$ operators $\mathfrak m_k$
in the context of $A_\infty$ algebra of Lagrangian submanifolds.
For a given compatible almost complex structure $J$, consider the moduli
space of stable maps of genus zero
$$
\CM_{k+1}(\beta;L) =\{ (w, (z_0,z_1, \cdots,z_k)) \mid
\overline \partial_J w = 0, \, z_i \in \partial D^2, \, [w] = \beta
\, \mbox{in }\, \pi_2(M,L) \}/\sim
$$
where $\sim$ is the conformal reparameterization of the disc $D^2$. The
expected dimension of this space is given by
\begin{equation}\label{eq:dim}
n+ \mu(\beta) - 3 + (k+1) = n+\mu(\beta) + k-2.
\end{equation}
Now given $k$ chains
$$
[P_1,f_1], \cdots,[P_k,f_k] \in C_*(L)
$$
of $L$ considered as \emph{currents} on $L$, we put the cohomological grading
$\mbox{deg} P_i = n - \dim P_i$ and consider the fiber product
$$
ev_0: \CM_{k+1}(\beta;L) \times_{(ev_1, \cdots, ev_k)}(P_1 \times
\cdots \times P_k) \to L.
$$
A simple calculation shows that the expected dimension of this chain is given by
$$
n + \mu(\beta) - 2 + \sum_{j=1}^k(\dim P_j + 1- n)
$$
or equivalently we have the expected degree
$$
\mbox{deg}\left[\CM_{k+1}(\beta;L) \times_{(ev_1, \cdots, ev_k)}(P_1
\times \cdots\times P_k),
ev_0\right] = \sum_{j=1}^n(\mbox{deg} P_j -1) + 2- \mu(\beta).
$$
For each given $\beta \in \pi_2(M,L)$ and $k = 0, \cdots$, we
define
$$
\mathfrak m_{k,\beta}(P_1, \cdots, P_k)
= \left[\CM_{k+1}(\beta;L) \times_{(ev_1, \cdots, ev_k)}(P_1 \times \cdots \times P_k),
ev_0\right]
$$
and $\mathfrak m_k = \sum_{\beta \in \pi_2(M,L)} \mathfrak m_{k,\beta}
\cdot q^\beta$
where $q^\beta = T^{\omega(\beta)} e^{\mu(\beta)/2}$ with $T,\,e$
formal parameters encoding the area and the Maslov index of $\beta$.
We provide $T$ with degree 0 and $e$ with $2$. Now we denote by $C[1]$
the completion of a \emph{suitably chosen countably generated} cochain
complex with $\Lambda_{0,nov}$ as
its coefficients with respect to the non-Archimedean topology. Then it follows that the map
$
\mathfrak m_k : C[1]^{\otimes k} \to C[1]
$
is well-defined, has degree 1 and continuous with respect to non-Archimedean topology.
We extend $\mathfrak m_k$ as a coderivation
$\widehat{\mathfrak m}_k: BC[1] \to BC[1]$
where $BC[1]$ is the completion of the direct sum $\oplus_{k=0}^\infty
B^kC[1]$ where $B^kC[1]$ itself is the completion of $C[1]^{\otimes k}$.
$BC[1]$ has a natural filtration defined similarly as \ref{eq:vlambda1}
Finally we take the sum
$$
\widehat d = \sum_{k=0}^\infty \widehat m_k : BC[1] \to BC[1].
$$
A main theorem then is the following coboundary property

\begin{thm}\label{algebra} Let $L$ be an arbitrary compact relatively
spin Lagrangian submanifold of an arbitrary tame symplectic manifold
$(M,\omega)$. The coderivation $\widehat d$ is a continuous map that
satisfies the $A_\infty$ relation $\widehat d \widehat d = 0$, and so
$(C,\frak m)$ is a filtered weak $A_\infty$ algebra.
\end{thm}
One might want to consider the homology of this huge complex but
if one naively takes homology of this complex itself, it will end up with
getting a trivial group, which is isomorphic to the ground ring
$\Lambda_{0, nov}$. This is because the $A_\infty$ algebra
associated to $L$ in \cite{FOOO} has the \emph{(homotopy) unit} :
if an $A_\infty$ algebra has
a unit, the homology of $\widehat d$ is isomorphic to its ground ring.

A more geometrically useful homology relevant to the Floer homology
is the $\mathfrak m_1$-homology (\ref{eq:m1cohom}) in this context,
which is the Bott-Morse version of the Floer cohomology for the pair
$(L,L)$. However in the presence of $\mathfrak m_0$,
$\widehat{\mathfrak m}_1 \widehat{\mathfrak m}_1 = 0$ no longer holds
in general. Motivated by Kontsevich's suggestion
\cite{konts:remark}, this led Fukaya-Oh-Ohta-Ono
to consider deforming Floer's original definition by a bounding
chain of the obstruction cycle arising from bubbling-off discs.
One can always deform the given (filtered) $A_\infty$ algebra $(C,\mathfrak m)$ by
an element $b \in C[1]^0$ by re-defining the $A_\infty$ operators as
$$
\mathfrak m_k^b(x_1,\cdots, x_k) = \mathfrak m(e^b,x_1, e^b,x_2, e^b,x_3, \cdots,
x_k,e^b)
$$
and taking the sum $\widehat d^b = \sum_{k=0}^\infty \widehat{\mathfrak m}_k^b$.
This defines a new weak $A_\infty$ algebra in general. Here we
simplify notations by writing
$$
e^b = 1 + b + b\otimes b + \cdots + b \otimes \cdots \otimes b +\cdots.
$$
Note that each summand in this infinite sum has degree 0 in $C[1]$ and
converges in the non-Archimedean topology if $b$ has positive
valuation, i.e., $v(b) > 0$.

\begin{prop} For the $A_\infty$ algebra $(C,\mathfrak m_k^b)$,
$\mathfrak m_0^b = 0$ if and only if $b$ satisfies
\begin{equation}\label{eq:MC}
\sum_{k=0}^\infty\frak m_k(b,\cdots, b) = 0.
\end{equation}
This equation is a version of \emph{Maurer-Cartan equation}
for the filtered $A_\infty$ algebra.
\end{prop}

\begin{defn}\label{boundchain}
Let $(C,\frak m)$ be a filtered weak $A_\infty$ algebra in general and
$BC[1]$ be its bar complex. An element $b \in C[1]^0 = C^1$ is called
a \emph{bounding cochain} if it satisfies the equation (\ref{eq:MC})
and $v(b) > 0$.
\end{defn}

In general a given $A_\infty$ algebra may or may not have a solution
to (\ref{eq:MC}).

\begin{defn}\label{unobstructed}
A filtered weak $A_\infty$ algebra is called \emph{unobstructed} if the equation
(\ref{eq:MC}) has a solution $b \in C[1]^0 = C^1$ with $v(b) > 0$.
\end{defn}
One can define a notion of homotopy equivalence between two
bounding cochains and et al as described in \cite{FOOO}. We denote by
$\CM(L)$ the set of equivalence classes of bounding cochains of $L$.

Once the $A_\infty$ algebra is attached to each Lagrangian submanifold $L$,
we then construct an \emph{$A_\infty$ bimodule} $C(L,L')$ for the pair
by considering operators
$$
\mathfrak n_{k_1,k_2}: C(L,L') \to C(L,L')
$$
defined similarly to $\mathfrak m_k$ : A typical generator of $C(L,L')$
has the form
$$
P_{1,1} \otimes \cdots, \otimes P_{1,k_1} \otimes [p,w] \otimes
P_{2,1} \otimes  \cdots \otimes P_{2,k_2}
$$
and then the image $\mathfrak n_{k_1,k_2}$ thereof is given by
$$
\sum_{[q,w']}\left[\left(\CM([p,w],[q,w']);P_{1,1},\cdots,P_{1,k_1};
P_{2,1},\cdots,P_{2,k_2}),ev_\infty\right)\right][q,w'].
$$
Here $\CM([p,w],[q,w']);P_{1,1},\cdots,P_{1,k_1};
P_{2,1},\cdots,P_{2,k_2})$ is the Floer moduli space
$$
\CM([p,w],[q,w']) = \bigcup_{[q,w'] = [q,w\#B]}\CM(p,q;B)
$$
cut-down by intersecting with the given chains $P_{1,i} \subset L$
and $P_{2,j} \subset L'$, and the evaluation map
$$
ev_\infty: \CM([p,w],[q,w']);P_{1,1},\cdots,P_{1,k_1};
P_{2,1},\cdots,P_{2,k_2}) \to \mbox{Crit}\CA
$$
is defined by $ev_\infty(u) = u(+\infty)$.
\begin{thm} Let $(L,L')$ be an arbitrary relatively spin pair
of compact Lagrangian submanifolds.
Then the family $\{\mathfrak n_{k_1,k_2}\}$ define
a left $(C(L),\mathfrak m)$ and right $(C(L'),\mathfrak m')$ filtered
$A_\infty$ bimodule structure on $C(L,L')$.
\end{thm}

In other words, each of the map $\mathfrak n_{k_1,k_2}$ extends to
a $A_\infty$ bimodule homomorphism $\widehat{\mathfrak n}_{k_1,k_2}$ and
if we take the sum
$$
\widehat d: = \sum_{k_1,k_2} \widehat{\mathfrak n}_{k_1,k_2}
: C(L,L') \to C(L,L'),
$$
$\widehat d$ satisfies the following coboundary property
\begin{prop}\label{bimodule} The map $\widehat d$ is a continuous map
and satisfies $\widehat d \widehat d = 0$.
\end{prop}

Again this complex is too big for the computational purpose
and we would like to consider the Floer homology by restricting the $A_\infty$
bimodule to a much smaller complex, an ordinary $\Lambda_{nov}$ module
$CF(L,L')$. However Floer's original definition again meets obstruction
coming from the obstructions cycles of either $L_0$, $L_1$ or of both.
We need to deform Floer's `boundary' map $\delta$
using suitable bounding cochains of $L, \, L'$. The bimodule $C(L,L')$
is introduced to perform this deformation coherently.
\par
In the case where both $L, \, L'$ are unobstructed, we can carry out this
deformation of $\mathfrak n$ by bounding chains $b_1 \in \CM(L)$ and $b_2 \in \CM(L')$
similarly as $\frak m^b$ above. Symbolically we can write the new operator
as
$$
\delta^{b_1,b_2}(x) = \mathfrak{\widehat n}(e^{b_1},x,e^{b_2}).
$$
\begin{thm} For each $b_1 \in \CM(L)$ and $b_2 \in \CM(L')$, the map
$\delta^{b_1,b_2}$ defines a continuous map
$
\delta^{b_1,b_2} : CF(L,L') \to CF(L,L')
$
that satisfies $\delta^{b_1,b_2}\delta^{b_1,b_2} = 0$.
\end{thm}

This theorem enables us to define the \emph{deformed Floer cohomology}

\begin{defn} For each $b \in \CM(L)$ and $b' \in \CM(L')$,
we define the \emph{$(b,b')$-Floer cohomology} of the pair $(L,L')$ by
$$
HF((L,b), (L',b');\Lambda_{nov}) = \frac{\ker \delta^{b_1,b_2}}
{\mbox{\rm im } \delta^{b_1,b_2}}.
$$
\end{defn}
\begin{thm}
The above cohomology remains isomorphic under the
Hamiltonian isotopy of $L, \, L'$ and under the
homotopy of bounding cochains $b, \, b'$.
\end{thm}

We refer to \cite{FOOO} and its revised version for all the details of
algebraic languages needed to make the statements in the above theorems
precise.

\subsection{Spectral sequence}

The idea of spectral sequence is quite simple to describe. One can follow
more or less the standard construction of spectral sequence on the filtered
complex e.g. in \cite{mccl}. One trouble to overcome in the
construction of spectral sequence on $(C(L),\delta)$ or $(C(L,L'),\delta)$
is that the general Novikov ring, in particular $\Lambda_{0,nov}$
is \emph{not} N\"oetherian and so the standard theorems on the N\"oetherian
modules cannot be applied. In addition, the Floer complex is
not bounded above which also makes the proof of convergence of the
spectral sequence somewhat tricky. We refer to \cite{FOOO} for complete discussions
on the construction of the spectral sequence and the study of their convergences.

However for the case of monotone Lagrangian submanifolds, the Novikov
ring becomes a field and the corresponding spectral sequence is
much simplified as originally carried out by Oh \cite{oh:imrn}
by a crude analysis of thick-thin decomposition of Floer moduli
spaces as two Lagrangian submanifolds collapse to one.
Then the geometric origin of the spectral sequence is the decomposition
of the Floer boundary map $\delta$ into
$\delta = \delta_0 + \delta_1 + \delta_2 + \cdots$
where each $\delta_i$ is the contribution coming from the Floer
trajectories of a given symplectic area in a way that the corresponding
area is increasing as $i \to \infty$. Here $\delta_0$ is the contribution
from the classical cohomology. In general this sequence
may not stop at a finite stage but it does for monotone Lagrangian
submanifolds. In this regard, we can roughly state the following general
theorem : \begin{itemize}
\item There exists a spectral sequence whose $E^2$-term is
isomorphic to the singular cohomology $H^*(L)$ and which
converges to the Floer cohomology $HF^*(L,L)$.
\end{itemize}
See \cite{oh:imrn} and \cite{FOOO} for the details of the monotone case and of
the general case respectively. The above decomposition also provides an algorithm to
utilize the spectral sequence in examples, especially when
the Floer cohomology is known as for the case of Lagrangian
submanifolds in $\C^n$. Here are some sample results.

\begin{thm}[Theorem II \cite{oh:imrn}]\label{ohimrnII} Let
$(M,\omega)$ be a tame symplectic manifold with $\dim M \geq 4$.
Let $L$ be a compact monotone
Lagrangian submanifold of $M$ and $\phi$ be a Hamiltonian
diffeomorphism of $(M,\omega)$ such that $L$ intersects $\phi(L)$
transversely. Then the followings hold :
\par
$1.$ If $\Sigma_L \geq n+2$,
$HF^k(L,\phi(L);\Z_2) \cong H^k(L;\Z_2)$ for all $k \mod
\Sigma_L$.
\par
$2.$ If $\Sigma_L = n+1$, the same is true for $k \neq 0, n \mod n+1$.
\end{thm}

\begin{thm}[Theorem III \cite{oh:imrn}]\label{ohimrnIII} Let
$L \subset \C^n$ be a compact monotone Lagrangian torus. Then we have
$\Sigma_L = 2$ provided $1 \leq  n \leq 24$.
\end{thm}

A similar consideration, using a more precise form of the spectral sequence from
\cite{FOOO}, proves

\begin{thm}\label{muneq0} Let $L\subset \C^n$ be a compact
Lagrangian embedding with $H^2(L;\Z_2) = 0$. Then
its Maslov class $\mu_L$ is not zero.
\end{thm}

The following theorem can be derived from Theorem E \cite{FOOO}
which should be useful for the study of intersection
properties of special Lagrangian submanifolds on Calabi-Yau manifold.

\begin{thm} Let $M$ be a Calabi-Yau manifold and $L$ be an
unobstructed Lagrangian submanifold with its Maslov class
$\mu_L = 0$ in $H^1(L;\Z)$. Then we have
$HF^i(L;\Lambda_{0,nov}) \neq 0$ for $i = 0, \, \dim L$.
\end{thm}

For example, any special Lagrangian
homology sphere satisfies all the hypotheses required in this theorem.
Using this result combined with some Morse theory argument,
Thomas and Yau \cite{thomas-yau} proved the following uniqueness result of
special Lagrangian homology sphere in its Hamiltonian isotopy class

\begin{thm}[Thomas-Yau] For any Hamiltonian isotopy class of
embedded Lagrangian submanifold $L$ with $H^*(L) \cong H^*(S^n)$,
there exists at most one smooth special Lagrangian representative.
\end{thm}

Biran \cite{biran} also used this spectral sequence
for the study of geometry of Lagrangian skeletons and
polarizations of K\"ahler manifolds.

\section{Displaceable Lagrangian submanifolds}

\begin{defn} We call a compact Lagrangian submanifold
$L \subset (M,\omega)$ displaceable if there exists a Hamiltonian
isotopy $\phi_H$ such that $L\cap \phi_H^1(L) = \emptyset$.
\end{defn}
One motivating question for studying such Lagrangian submanifolds
is the following well-known folklore conjecture in symplectic geometry.
\begin{conj}[Maslov Class Conjecture] \label{maslov}
Any compact Lagrangian embedding in
$\C^n$ has non-zero Maslov class.
\end{conj}
Polterovich \cite{pol} proved the conjecture in dimension $n=2$ whose proof uses a loop $\gamma$
realized by the boundary of Gromov's holomorphic disc
constructed in \cite{gromov}. Viterbo proved
this conjecture for any Lagrangian torus in $\C^n$ by a different
method using the critical point theory on the free loop spaces
of $\C^n$ \cite{viterbo1}. Also see Theorem \ref{muneq0} in the previous section
for $L$ with $H^2(L;\Z_2)\neq 0$.

It follows from definition that $HF^*(L,\phi_H^1(L))= 0$
for a displaceable Lagrangian submanifold $L$ \emph{whenever
$HF^*(L,\phi_H^1(L))$ is defined}.
An obvious class of displaceable Lagrangian submanifolds are those
in $\C^n$. This simple observation, when combined with
the spectral sequence described in the previous section, provides
many interesting consequences on the symplectic topology
of such Lagrangian submanifolds as illustrated by
Theorem \ref{ohimrnII} and \ref{ohimrnIII}.

Some further amplification of this line of
reasoning was made by Biran
and Cieliebak \cite{biran-ciel1} for the study of topology
of Lagrangian submanifolds in (complete) \emph{sub-critical Stein manifolds}
$(V,J)$ or a symplectic manifold $M$ with such $V$ as a factor.
They cooked up some class of Lagrangian submanifolds in such symplectic
manifolds with suitable condition on the first Chern class of $M$ under which
the Lagrangian submanifolds become monotone and satisfy
the hypotheses in Theorem \ref{ohimrnII}. Then applying this theorem, they
derived restrictions on the topology of such Lagrangian submanifolds,
e.g., some \emph{cohomological
sphericality} of such Lagrangian submanifolds
(see Theorem 1.1 \cite{biran-ciel1}).

Recently Fukaya \cite{FuLag} gave a new construction of the
$A_\infty$-structure described in the previous section as a deformation
of the differential graded algebra of de Rham complex of $L$
associated to a natural solution to the Maurer-Cartan
equation of the Batalin-Vilkovisky structure discovered by
Chas and Sullivan \cite{chas-sull} on the loop space. In this way,
Fukaya combined  Gromov and Polterovich's pseudo-holomorphic curve
approach and Viterbo's loop space approach \cite{viterbo1}, and proved several new
results on the structure of Lagrangian embeddings in $\C^n$.
The following are some sample results proven by this method \cite{FuLag} :
\par
1. If $L$ is spin and aspherical in $\C^n$ then a finite cover $\widetilde L$ of $L$ is
homotopy equivalent to a product $S^1 \times \widetilde L'$.
Moreover the Maslov index of $[S^1] \times [point]$ is $2$.
\par
2. If $S^1 \times S^{2n}$ is embedded as a Lagrangian
submanifold of $\C^{2n+1}$, then the Maslov
index of $[S^1]\times [point]$ is 2.
\par
There is also the symplectic field theory approach to the proof
of the result 1 above for the case of torus $L = T^n$ as
Eliashberg explained to the authors \cite{eliash2}.
Eliashberg's scheme has been further detailed by Cieliebak
and Mohnke \cite{ciel-mohn}.
The result 1 for $T^n$ answers affirmatively to Audin's
question \cite{audin} on the minimal Maslov number of the
embedded Lagrangian torus in $\C^n$ for general $n$. Previously
this was known only for $n = 2$ \cite{pol}, \cite{viterbo1} and
for monotone Lagrangian tori \cite{oh:imrn} (see Theorem \ref{ohimrnIII}).

\section{Applications to mirror symmetry}

Mirror symmetry discovered in the super string theory
attracted much attention from many (algebraic) geometers since it
made a remarkable prediction on the relation between
the number of rational curves on a Calabi-Yau 3-fold $M$ and
the deformation theory of complex structures of another
Calabi-Yau manifold $M^{\dag}$.

\subsection{Homological mirror symmetry}

Based on Fukaya's construction of the $A_\infty$ category of symplectic
manifolds \cite{Ainfold},
Kontsevich \cite{konts:Zur} proposed a conjecture on the relation between
the category ${\rm Fuk}(M)$ of $(M,\omega)$ and the derived category of coherent
sheaves ${\rm Coh}(M^{\dag})$ of $M^{\dag}$, and extended the
mirror conjecture in a more conceptual way.
This extended version is called the \emph{homological mirror symmetry},
which is closely related to the D-brane duality studied
 much in physics. Due to the obstruction phenomenon
we described in \S 3.4, the original construction in \cite{Ainfold}
requires some clarification of the definition of ${\rm Fuk}(M)$.
The necessary modification has been completed in \cite{FOOO,MsurvII}.

For the rest of this subsection, we will formulate a precise mathematical
conjecture of homological mirror symmetry.
Let $(M,\omega)$ be an integral symplectic manifold i.e.,
one with $[\omega] \in H^2(M;\Z)$. For such $(M,\omega)$,
we consider a family of complexified symplectic structures
$\quad
M_{\tau} = (M, -\sqrt{-1}\tau \omega)\quad
$
parameterized by $\tau \in \mathfrak h$ where $\mathfrak h$
is the upper half plane. The mirror of this family is expected
to be a family of complex manifolds $M^{\dag}_{q}$ parameterized by
$q = e^{\sqrt{-1}\tau}\in D^2\setminus \{0\}$, the punctured disc.
Suitably `formalizing' this family at $0$, we obtain
a scheme $\mathfrak M^{\dag}$ defined over the ring
$\Q[[q]][q^{-1}]$. We identify $\Q[[q]][q^{-1}]$ with a sub-ring
of the universal Novikov ring $\Lambda_{nov}$ defined in
subsection \ref{subsec:off-shell}.
The ext group $\text{Ext}({\mathcal E}_0,{\mathcal E}_1)$ between the
coherent sheaves ${\mathcal E}_i$ on $\mathfrak M^{\dag}$ is a module over
$\Q[[q]][q^{-1}]$.

We consider the quadruple ${\mathcal L} = (L,s,d,[b])$, which we call
a \emph{Lagrangian brane}, that satisfies the following data :
\par
1. $L$ a Lagrangian submanifold of $M$
such that the Maslov index of $L$ is zero and
$[\omega] \in H^2(M,L;\Z)$. We also enhance $L$ with flat complex
line bundle on it.
\par
2. $s$ is a spin structure of $L$.\par
3. $d$ is a grading in the sense of \cite{konts:Zur}, \cite{seidel1}.
\par
4. $[b] \in \mathcal M(L)$ is a bounding cochain described in subsection
\ref{subsec:obstruction}.

\begin{conj}
To each Lagrangian brane $\mathcal L$ as above, we can associate
an object ${\mathcal E}(\mathcal L)$ of the derived category of coherent sheaves on
the scheme $\mathfrak M^{\dag}$ so that the following holds :
\par
$1.$ There exists a canonical isomorphism.
$$
HF({\mathcal L}_1,{\mathcal L}_2) \cong \text{\rm Ext}({\mathcal E}({\mathcal L}_1),{\mathcal E}
({\mathcal L}_2))\otimes_{\Q[[q]][q^{-1}]}\Lambda_{nov}
$$
\par
$2.$ The isomorphism in 1 is functorial :
Namely the product of Floer cohomology is mapped to the
Yoneda product of Ext group by the isomorphism 1.
\end{conj}

The correct Floer cohomology $HF({\mathcal L}_1,{\mathcal L}_2)$
used in this formulation of the conjecture
is given in \cite{FOOO} (see \S 3.4 for a brief description).
The spin structure in ${\mathcal L}$ is needed to define
orientations on the various moduli spaces involved
in the definition of Floer cohomology and
the grading $d$ is used to define an absolute integer grading
on $HF({\mathcal L}_1,{\mathcal L}_2)$. We refer readers to \cite{FOOO} \S 1.4,
\cite{abel} for the details of construction and for more references.
\par
We now provide some evidences for this conjecture.
A conjecture of this kind was first
observed by Kontsevich in \cite{konts:Zur} for the case of an elliptic
curve $M$, which is further explored by Polischchuk-Zaslow \cite{PZ},
and by Fukaya in \cite{abel} for the case
when $M$ is a torus (and so $M^{\dag}$ is also a torus)
and $L \subset M$ is an affine sub-torus. In fact, in these cases
one can use the convergent power series for the
formal power series or the Novikov ring.
Kontsevich-Soibelman \cite{KS} gave an alternative proof,
based on the adiabatic degeneration result of the authors \cite{foh},
for the case where $L$ is an etal\'e covering of the base torus
of the Lagrangian torus fibration $M = T^{2n} \to T^n$.
Seidel  proved Conjecture 5.1 for the quartic surface $M$ \cite{seidel2}.
\par

\subsection{Toric Fano and Landau-Ginzburg
correspondence}

So far we have discussed the case of Calabi-Yau
manifolds (or a symplectic manifold $(M,\omega)$ with $c_1(M) = 0$).
The other important case
that physicists studied much is the case of toric Fano
manifolds, which physicists call the correspondence between
the $\sigma$-model and the Landau-Ginzburg model.
Referring readers to \cite{hori} and \cite{hori-vafa} for detailed
physical description of this correspondence, we briefly describe
an application of machinery developed in \cite{FOOO} for
an explicit calculation of Floer cohomology of Lagrangian
torus orbits of toric Fano manifolds.  We will focus on
the correspondence of the $A$-model of a toric Fano manifold and the
$B$-model of Landau-Ginzburg model of its mirror.
We refer to \cite{HIV} for the other side of the correspondence
between the toric Fano $B$-model and the Landau-Ginzburg $A$-model.

According to \cite{FOOO}, the obstruction cycles
of the filtered $A_\infty$ algebra associated to a Lagrangian
submanifold is closely related to $\mathfrak m_0$. This $\mathfrak m_0$,
by definition, is defined by a collection of the (co)chains
$
[\CM_1(\beta), ev_0]
$
for all $\beta \in \pi_2(M,L)$. More precisely, we have
\begin{equation}\label{eq:m0}
\mathfrak m_0(1) = \sum_{\beta \in \pi_2(M,L)}[\CM_1(\beta), ev_0]
\cdot T^{\omega(\beta)} q^{\mu(\beta)/2} \in C^*(L) \otimes
\Lambda_{0,nov}.
\end{equation}
This is the sum of all genus zero instanton contributions with one marked
point.

On the other hand, based on a $B$-model calculations,
Hori \cite{hori}, Hori-Vafa \cite{hori-vafa}
proposed some correspondence between the instanton contributions
of the $A$-model of toric Fano manifolds and the Landau-Ginzburg
potential of the $B$-model of its mirror. This correspondence
was made precise by Cho and Oh \cite{cho-oh}. A  description of
this correspondence is now in order.

First they proved the following
\begin{thm}
$[\CM_1(\beta),ev_0]=[L]$ as a chain, for every $\beta \in \pi_2(M,L)
$ with $\mu(\beta) = 2$ and so
$\frak m_0(1) = \lambda [L]$ for some $\lambda \in \Lambda_{0,nov}$.
\end{thm}
It had been previously observed in Addenda of \cite{oh:cpam} for the monotone case
that Floer cohomology $HF^*(L,L)$ is defined even when the
minimal Maslov number $\Sigma_L = 2$. Using the same argument
Cho and Oh proved that $HF^*(L, L;\Lambda_{0,nov})$
is well-defined for the torus fibers of toric Fano manifolds
\emph{without deforming} Floer's `boundary' map, at least for
the \emph{convex} case. We believe this convexity condition can be removed.
More specifically,
they proved $\frak m_1\frak m_1 = 0$.
This is because in (\ref{eq:m0m1}), the last two
terms cancel each other \emph{if $\frak m_0(1) = \lambda e$ is a
multiple of the unit $e = [L]$} and then (\ref{eq:m1m0=0}) implies that
$\frak m_0(1)$ is a $\frak m_1$-cycle.
We refer readers to the revision of \cite{FOOO}
for more discussion on this case, in which
any filtered $A_\infty$ algebra deformable to such a case
is called \emph{weakly unobstructed}.

In fact, Cho and Oh obtained the explicit formula
\begin{equation}\label{eq:o}
\mathfrak m_0(1) = \sum_{i=1}^N h^{v_j}e^{\sqrt{-1}\langle \nu, v_j\rangle}
T^{\omega(\beta_j)}[L]\cdot q
\end{equation}
after including flat line bundles attached to $L$ and computing
precise formulae for the area $\omega(\beta_j)$'s.
Here $h^{v_j} = e^{\sqrt{-1}\langle \nu, v_j\rangle}$ is the holonomy
of the flat line bundle and $\omega(\beta_j)$ was calculated
explicitly in \cite{cho-oh}. Denote by $\nu = (\nu_1, \cdots, \nu_N)$ the
holonomy vector of the line bundle appearing in the description of
linear $\sigma$-model \cite{hori-vafa}.

On the other hand, the Landau-Ginzburg potential is given by the formula
$$
\sum_{i=1}^N \exp(-y_i - \langle \Theta, v_i \rangle) =: W(\Theta)
$$
for the mirror of the given toric manifold (see \cite{hori-vafa} for
example).

\begin{thm}[\cite{cho-oh}] Let $A \in \mathfrak t^*$ and denote
$\Theta = A - \sqrt{-1}\nu$. We denote
by $\mathfrak m^\Theta$ the $A_\infty$ operators associated to the torus fiber
$T_A = \pi^{-1}(A)$ coupled with the flat line bundle whose holonomy vector is
given by $\nu \in (S^1)^N$ over the toric fibration $\pi: X \to \mathfrak t^*$.
Under the substitution of $T^{2\pi} = e^{-1}$ and ignoring
the harmless grading parameter $q$, we have the exact correspondence
\begin{eqnarray}
\mathfrak m_0^\Theta & \longleftrightarrow & W(\Theta) \\
\mathfrak m_1^\Theta (pt) & \longleftrightarrow & dW(\Theta) = \sum_{j=1}^n \frac{\partial
W}{\partial \Theta_j}(\Theta) d\Theta_j
\end{eqnarray}
under the mirror map given in \emph{\cite{hori-vafa}}.
\end{thm}

Combined with a theorem from \cite{cho-oh}
which states that $HF^*(\CL,\CL) \cong H^*(L;\C)\otimes \Lambda_{0,nov}$
whenever $\mathfrak m_1^\Theta (pt) = 0$, this theorem confirms the
prediction made by Hori \cite{hori}, Hori-Vafa \cite{hori-vafa} about
the Floer cohomology of Lagrangian torus fibers.
This theorem has been
further enhanced  by Cho \cite{cho} who relates the higher
derivatives of $W$ with the higher Massey products $\mathfrak m_k^\Theta$.
For example, Cho proved that the natural product structure on
$HF^*(L,L)$ is not isomorphic to the cohomology ring $H^*(T^n)\otimes \Lambda_{0,nov}$
but isomorphic to the Clifford algebra associated to the quadratic
form given by the Hessian of the potential $W$ under the mirror
map. This was also predicted by physicists (see \cite{hori}).

\vskip0.2in



\begin{thebibliography}{FOOO}

\bibitem[Au]{audin} Audin, M., Fibres normaux d'immersions en dimension
double, points doubles d'immersions Lagrangiennes et plongements totalement
r\'eeles, \emph{Comment. Math. Helv.} \textbf{63} (1988), 593--623.

\bibitem[Ba]{banyaga} Banyaga, A., Sur la structure du groupe des
diff\'eomorphismes qui pr\'eservent une forme symplectique,
\emph{Comm. Math. Helv.} \textbf{53}  (1978),  174--227.

\bibitem[Bi]{biran} Biran, P., Lagrangian non-intersections,
preprint 2004, arXiv:math.SG/0412110.

\bibitem[BC]{biran-ciel1} Biran, P., Cieliebak, K.,
Lagrangian embeddings into subcritical Stein manifolds,
\emph{Isreal J. Math.} \textbf{127} (2002), 221--244.

\bibitem[CS]{chas-sull} Chas, M., Sullivan D., String topology,
preprint 1999, math.GT/9911159.

\bibitem[Cho]{cho} Cho, C.-H., Products of Floer cohomology of
Lagrangian torus fibers in toric Fano manifolds,
Comm. Math. Phys.  (to appear).

\bibitem[CO]{cho-oh} Cho, C.-H., Oh, Y.-G. Floer cohomology and disc
instantons of Lagrangian torus fibers in toric Fano manifolds,
submitted, arXiv:math.SG/0308225.

\bibitem[CM]{ciel-mohn} Cieliebak, K., Mohnke, K. a talk by
Cieliebak in Banach Center, Warsaw, 2004.

\bibitem[EH]{ekel-hofer} Ekeland, I., Hofer, H., Symplectic topology
and Hamiltonian dynamics \& II, \emph{Math. Z.} \textbf{200} (1990), 355--378,
\& \textbf{203} (1990), 553-567.

\bibitem[El1]{eliash1} Eliashberg, Y., A theorem on the structure
of wave fronts and applications in symplectic topology,
\emph{Funct. Anal. and its Appl.} \textbf{21} (1987), 227--232.

\bibitem[El2]{eliash2} Eliashberg, Y., private communication, 2001.

\bibitem[EnP]{entov-pol1}  Entov, M., Polterovich, L., Calabi
quasimorphism and quantum homology, \emph{Internat. Math. Res. Notices} no
\textbf{30} (2003),  1635--1676.

\bibitem[Fa]{fathi} Fathi, A., Structure of the group of
homeomorphisms preserving a good measure on a compact manifold,
\emph{Ann. Scient. \`Ec. Norm. Sup.} \textbf{13} (1980), 45--93.

\bibitem[Fl1]{floer:intersect} Floer, A., Morse theory for Lagrangian
interesections, \emph{J. Differ. Geom.} \textbf{28} (1988), 513--547.

\bibitem[Fl2]{floer:fixed} Floer, A., Symplectic fixed points and
holomorphic spheres, \emph{Commun. Math. Phys.}  \textbf{120} (1989),
575--611.

\bibitem[Fu1]{Ainfold} Fukaya, K.,
Morse homotopy, $A\sp \infty$-category, and Floer homologies,
\emph{Proceedings of
GARC Workshop on Geometry and Topology} '93 , 1--102, Lecture
Notes Ser., \textbf{18}, Seoul Nat. Univ., Seoul, (1993).

\bibitem[Fu2]{MsurvII} Fukaya, K.,Floer homology and mirror symmetry. II,
31--127, \emph{Adv. Stud. Pure Math., 34, Math. Soc. Japan,
Tokyo}, (2002).

\bibitem[Fu3]{abel} Fukaya, K., Mirror symmetry of abelian varieties and
multi-theta functions, \emph{J. Algebraic Geom.}  \textbf{11} (2002), 393--512.

\bibitem[Fu4]{FuLag} Fukaya, K., Application of Floer Homology of Langrangian
Submanifolds to Symplectic Topology, \emph{Morse Theoretic Methods in
Nonlinear Analysis and in Symplectic Topology} (2006).

\bibitem[FOh]{foh} Fukaya, K., Oh, Y.-G., Zero-loop open
strings in the cotangent bundle and Morse homotopy,  \emph{Asian J.
Math.}  \textbf{1} (1997),  96--180.

\bibitem[FOOO]{FOOO} Fukaya, K., Oh. Y.-G., Ohta, H., Ono, K.
\emph{Lagrangian intersection Floer homology - anomaly
and obstruction}, preprint 2000, http://www.math.kyoto-u.ac.jp/~fukaya/
; a revision in preparation.

\bibitem[FOn]{fon} Fukaya, K., Ono, K.,  Arnold conjecture and
Gromov-Witten invariants, \emph{Topology}  \textbf{38}  (1999),
933--1048.

\bibitem[GJ]{GJ} Getzler, E., Jones, D.S., $A_\infty$ algebra and
cyclic bar complex, \emph{Illinois J. Math.} \textbf{34} (1990), 256--283.

\bibitem[Gr]{gromov}  Gromov, M., Pseudo-holomorphic curves in
symplectic manifolds,  \emph{Invent. Math.}  \textbf{82} (1985), 307--347.

\bibitem[H]{hofer} Hofer, H. On the topological properties of
symplectic maps, \emph{Proc. Royal Soc. Edinburgh}  \textbf{115} (1990),
25--38.

\bibitem[HS]{hofer-sal}  Hofer, H., Salamon, D., Floer homology and
Novikov rings, in the \emph{Floer Memorial Volume}, Hofer, H. et
al, Birkha\"user, 1995, pp 483--524.

\bibitem[Ho]{hori} K. Hori, Linear models in supersymmetric D-branes,
Proceedings of the KIAS conference, \emph{Mirror Symmetry and Symplectic Geometry},
111--186, (Seoul 2000), eds by K. Fukaya, Y.-G. Oh, K. Ono and G. Tian,
World Sci. Publishing, River Edge, New Jersey, 2001.

\bibitem[HIV]{HIV} Hori, K., Iqbal, A., Vafa, C., D-branes and
mirror symmetry, preprint, hepth/0005247.

\bibitem[HV]{hori-vafa}
K. Hori and C. Vafa, Mirror symmetry, preprint, 2000,
hep-th/0002222.

\bibitem[K1]{konts:Zur} Kontsevich, M., Homological algebra of mirror symmetry.
\emph{Proceedings of the International Congress of Mathematicians, Vol. 1,
2 (Z\"uich, 1994), 120--139, Birkh\"aser, Basel, } 1995.

\bibitem[K2]{konts:remark} Kontsevich, M., private communication, 1997.

\bibitem[KS]{KS} Kontsevich, M., Soibelman. Y., Homological mirror symmetry
and torus fibrations. \emph{Symplectic geometry and Mirror Symmetry},
(Seoul 2000), 203--263, World Sci. Publishing, River Edge, New Jersey, 2001.

\bibitem[LT1]{liu-tian1} Liu, G., Tian, G., Floer homology
and Arnold conjecture, \emph{J. Differ. Geom.} \textbf{49} (1998), 1--74.

\bibitem[Mc]{mccl} McCleary, J., \emph{User's Guide to Spectral Sequences},
Math. Lec. Series \textbf{12}, Pulish or Perish, Wilmington, 1985.

\bibitem[N]{novikov1} Novikov, S.P., Multivalued functions and
functionals. An analogue of the Morse theory, (Russian) \emph{Dokl.
Akad. Nauk SSSR} \textbf{260} (1981), no. 1, 31--35.

\bibitem[Oh1]{oh:cpam} Oh, Y.-G., Floer cohomology of Lagrangian
intersections and pseudo-holomorphic discs, I \& II, \emph{Comm. Pure
and Applied Math.} \textbf{46} (1993), 949 -- 994 \& 995 -- 1012 ;
Addenda, ibid, \textbf{48} (1995), 1299 - 1302

\bibitem[Oh2]{oh:newton} Oh, Y.-G., Relative Floer and quantum cohomology and the symplectic
topology of Lagrangian submanifolds,  \emph{Contact and Symplectic Geometry},
 eds. by C. B. Thomas,  201--267,
Cambridge University Press, Cambridge, England, 1996.

\bibitem[Oh3]{oh:imrn}  Oh, Y.-G., Floer cohomology, spectral
sequences, and the Maslov class of Lagrangian embeddings,
\emph{Internat. Math. Res. Notices}, no 7 (1996),  305--346.

\bibitem[Oh4]{oh:jdg}  Oh, Y.-G., Symplectic topology as the
geometry of action functional, I, \emph{J. Differ. Geom.}  \textbf{46} (1997),
499--577.

\bibitem[Oh5]{oh:cag}  Oh, Y.-G., Symplectic topology as the
geometry of action functional, II,  \emph{Commun. Anal. Geom.}  \textbf{7} (1999),
1--55.

\bibitem[Oh6]{oh:ajm1} Oh, Y.-G., Chain level Floer theory and
Hofer's geometry of the Hamiltonian diffeomorphism group,
\emph{Asian J. Math.}  \textbf{6}  (2002),  579--624 ;
{\em Erratum}  \textbf{7} (2003), 447-448.

\bibitem[Oh7]{oh:ajm2}  Oh, Y.-G., Spectral invariants and length
minimizing property of Hamiltonian paths,  \emph{Asian J. Math.}
\textbf{9} (2005), 1--18.

\bibitem[Oh8]{oh:alan}  Oh, Y.-G., Construction of spectral
invariants of Hamiltonian paths on closed symplectic manifolds,
in \emph{The Breadth of Symplectic and Poisson Geometry},
Prog. Math. \textbf{232}, 525--570, Birkh\"auser, Boston, 2005.

\bibitem[Oh9]{oh:dmj} Oh, Y.-G., Spectral invariants, analysis
of the Floer moduli space and geometry of Hamiltonian
diffeomorphisms, \emph{Duke Math. J.} (to appear), math.SG/0403083.

\bibitem[Oh10]{oh:minimax}  Oh, Y.-G., Floer mini-max theory, the
Cerf diagram, and the spectral invariants, 2004 preprint,
math.SG/0406449.

\bibitem[Oh11]{oh:smoothing} Oh, Y.-G., $C^0$-coerciveness of Moser's
problem and smoothing of area preserving diffeomorphisms, submitted,
2005.

\bibitem[Oh12]{oh:hameo2} Oh, Y.-G., The group of Hamiltonian
homeomorphisms and topological Hamiltonian flows, preprint, December 2005.

\bibitem[OM]{oh:hameo1} Oh, Y.-G., M\"uller, S., The group of Hamiltonian
homeomorphisms and $C^0$-symplectic topology, submitted, a revision, December
2005, math.SG/0402210.

\bibitem[Ot]{ohta} Ohta, H., Obstruction to and deformation of
Lagrangian intersection Floer cohomology,
\emph{Mirror Symmetry and Symplectic Geometry}, (Seoul 2000), 281–-309,
World Sci. Publishing, River Edge, New Jersey, 2001.

\bibitem[On]{ono} Ono, K., On the Arnold conjecture for weakly
monotone symplectic manifolds, \emph{Invent. Math.} \textbf{119} (1995),
519--537.


\bibitem[P]{pol} Polterovich, L., The Maslov class of Lagrange
surfaces and Gromov's pseudo-holomorphic curves, \emph{Trans. Amer.
Soc.} \textbf{325} (1991), 241--248

\bibitem[PZ]{PZ} Polishchuk, A and Zaslow. E.,
Categorical mirror symmetry: the elliptic curve, \emph{Adv. Theor. Math. Phys.}
\textbf{2} (1998), no. 2, 443--470.

\bibitem[Ru]{ruan} Ruan, Y., Virtual neighborhood and
pseudo-holomorphic curves, \emph{Turkish J. Math.} \textbf{23} (1999),
161--231

\bibitem[Sc]{schwarz}  Schwarz, M., On the action spectrum for closed
symplectically aspherical manifolds, \emph{Pacific J. Math.} \textbf{193} (2000),
419--461

\bibitem[Se1]{seidel1} Seidel, P., Graded Lagrangian
submanifolds, \emph{Bull. Soc. Math. France} \textbf{128} (2000), 103--149.

\bibitem[Se2]{seidel2} Seidel, P., Homological mirror symmetry for the quartic
surface, preprint, 2003,  math.SG/0310414.

\bibitem[Si]{silva} De Silva, V. Products on symplectic Floer homology,
Thesis, Oxford Universiy, 1997.

\bibitem[St]{stasheff} Stasheff, J., Homotopy associativity of H-Spaces I \& II,
\emph{Trans. Amer. Math. Soc.} \textbf{108} (1963), 275--312 \& 293--312

\bibitem[TY]{thomas-yau} Thomas, R., S.-T. Yau, Special Lagrangians, stable
bundles and mean curvature flow, \emph{Comm. Anal. Geom.} \textbf{10} (2002),  no. 5,
1075--1113.

\bibitem[V1]{viterbo1} Viterbo, C., A new obstruction to embedding Lagrangian
tori, \emph{Invent. Math.} \textbf{100} (1990), 301--320.

\bibitem[V2]{viterbo2}  Viterbo, C., Symplectic topology as the
geometry of generating functions,  \emph{Math. Ann.}  \textbf{292}
(1992),  685--710.

\bibitem[V3]{viterbo3} Viterbo, C., On the uniqueness of generating
Hamiltonian for continuous limits of Hamiltonians flows, preprint, 2005,
math.SG/0509179

\end{thebibliography}
\end{document}